\documentclass[12pt]{article}
\usepackage{times}
\usepackage{amscd}
\usepackage{amssymb}

\newtheorem{theorem}{Theorem}[section]
\newtheorem{lemma}{Lemma}[section]

\def\ad{\mathop{\mathrm {ad}}\nolimits}
\def\Aut{\mathop{\mathrm {Aut}}\nolimits}
\def\Ad{\mathop{\mathrm {Ad}}\nolimits}
\def\End{\mathop{\mathrm {End}}\nolimits}
\def\Hom{\mathop{\mathrm {Hom}}\nolimits}
\def\Der{\mathop{\mathrm {Der}}\nolimits}
\def\Coh{\mathop{\mathrm {Coh}}\nolimits}
\def\cent{\mathop{\mathrm {cent}}\nolimits}
\def\GL{\mathop{\mathrm {GL}}\nolimits}
\def\Res{\mathop{\mathrm {Res}}\nolimits}

\begin{document}

\pagestyle{headings}
\markboth{Do Ngoc Diep}{Discrete Series for Loop Groups}

\title{{\sc Discrete Series for Loop Groups}.I\\
An Algebraic Realization of\\  Standard Modules
\footnote{Date: September 20, 1998, A revised version of {\bf
92-015}(1992), pp IV.1-IV.16,
Sonderforschungbereich 343, Universit\"at Bielefeld}}
\author{\sc Do Ngoc Diep \footnote{The author was supported in part by
the Alexander von Humboldt Foundation, the Abdus Salam ICTP and the Vietnam
National Foundation for Reserach in Fundamental Sciences}}
\date{{\bf Institute of Mathematics}\\ National Center for Science and Technology\\
P. O. Box 631, Bo Ho\\ VN-10.000  Hanoi\\
Vietnam  \footnote{Forwarded Email: {\tt
dndiep@member.ams.org}} \\
Email: {\tt dndiep@ioit.ncst.ac.vn }}
\maketitle

\begin{abstract} \noindent In this paper we consider the category
${\mathcal C} (\tilde{\mathfrak k}, \tilde{H})$ of the $(\tilde{\mathfrak
k},{\tilde H})$-modules, including all the Verma modules, where
${\mathfrak k}$ is some compact Lie algebra and H some Cartan subgroup,
$\tilde{\mathfrak k}$ and ${\tilde H}$ are the corresponding affine Lie
algebra and the affine Cartan group, respectively. To this category we
apply the Zuckermann functor and its derivatives. By using the determinant
bundle structure, we prove the natural duality of the Zuckermann derived
functors, deduce a Borel-Weil-Bott type theorem on decomposition of the
nilpotent part cohomology.  \end{abstract}

\section{Introduction}

Let us consider loop groups associated with compact Lie groups. 
Considering the (infinite dimensional) affine analogy of the classical
theories for finite dimensional Lie groups [1]-[3], we have in our
situation the following diagram of categories and functors

$$\begin{picture}(100,200)(10,20)
\put(-90,120){${\mathcal C}(\tilde{\mathfrak
h},\tilde{H})\quad\longrightarrow
\quad{\mathcal C}
(\tilde{\mathfrak b},\tilde{H})\quad{\buildrel H \over \longrightarrow}\quad 
{\mathcal O}\quad\cong \quad \oplus_{\chi} {\mathcal O}_{\chi}$}

\put(60,200){${\mathcal C}(\tilde{\mathfrak k},\tilde{H})\quad{\buildrel
R^iS
\over \longrightarrow} \quad {\mathcal C}
(\tilde{\mathfrak k},\tilde{H})^\approx\circlearrowright$}
\put(-60,140){\vector(2,1){100}}
\put(-10,180){I}
\put(80,160){$\bigcup$}	
\put(100,90){${\mathcal A}
(\tilde{\mathfrak k},\tilde{B}) = {\mathcal K}^\sim\circlearrowright $}
\put(110,60){$\cup$}
\put(80,30){$\{ M_{\lambda}= U(\tilde{\mathfrak k}) \otimes_{U(\tilde
{\mathfrak b})} L_{\lambda - \tilde{\rho}} \}$}
\put(100,15){(Verma \enskip modules)}
\end{picture}$$

where 

$I(W) := \Hom_{U(\tilde{\mathfrak h})}(U(\tilde{\mathfrak k}),W)[\tilde{\mathfrak
h}],$

$H(W) := \Hom_{U(\tilde{\mathfrak b})}(U(\tilde{\mathfrak k}),W)[\tilde{\mathfrak
h}],$

$S(V)  := V[\tilde{\mathfrak k}]$, the Zuckerman functor,

$V^\sim$ is the maximal $U(\tilde{\mathfrak h})$-locally finite and
$\tilde{H}$-semisimple module in the algebraic dual $V^*$,

$V^\approx$ the maximal $U(\tilde{\mathfrak k})$-locally finite and $\tilde{\mathfrak
k}$-semisimple module in $V^*$,

${\mathcal C}(\tilde{\mathfrak h},\tilde{H})$ the category of finite dimensional
semisimple $(\tilde{\mathfrak h},\tilde{H})$-modules,

${\mathcal C} (\tilde{\mathfrak b},\tilde{H})$ the category of the
$U(\tilde{\mathfrak h})$-locally finite, $\tilde{H}$-semisimple and
$U(\tilde{\mathfrak n}_+)$-locally finite $(\tilde{\mathfrak b},\tilde{H})$-modules,

${\mathcal C} (\tilde{\mathfrak k},\tilde{H})$ the category of the
$U(\tilde{\mathfrak h})$-locally finite, $\tilde{H}$-semisimple $(\tilde{\mathfrak
k},\tilde{H})$-modules,

${\mathcal O}$ the category of $U(\tilde{\mathfrak h})$-locally finite,
$\tilde{H}$-semisimple $U(\tilde{\mathfrak n}_{\pm})$-locally finite
$(\tilde{\mathfrak k},\tilde{H})$-modules,

${\mathcal A} (\tilde{\mathfrak k},\tilde{B})$ the category of $U(\tilde{\mathfrak
h})$-locally finite, $\tilde{H}$-semisimple finite dimensional
$\tilde{H}$-isotropic and $U(\tilde{\mathfrak n}_{\pm})$-locally finite
$(\tilde{\mathfrak k},\tilde{B})$-modules,

${\mathcal C} (\tilde{\mathfrak k},\tilde{K})$ the category of $(\tilde{\mathfrak
k},\tilde{K})$-modules. 

The functor $I$ shows that {\sl the category ${\mathcal C}
(\tilde{\mathfrak k},
\tilde{H})$ has enough injective objects}; hence, one can define the
Zuckerman
derived functors $$R^iS(V) := \bigoplus_{\lambda \in \tilde{P}} H^i(\tilde{\mathfrak
k},\tilde{H}; 
  V_{\lambda}^* \otimes V) \otimes V_{\lambda}.$$ For the objects from the
subcategory ${\mathcal O}$ \enskip one can restrict on the injective in
${\mathcal O}$ \enskip resolutions, i.e. $$R^i(S|_{\mathcal O})  \cong
(R^iS)|_{\mathcal O}.$$ By using the determinant bundle see for example
\cite{ps} structure $$\det (\tilde{\mathfrak k}/\tilde{\mathfrak h}) =
\wedge^{max}(\tilde{\mathfrak k}/\tilde{\mathfrak h}),$$ we can easily
prove the natural duality for relative Lie cohomology
$$H^i(\tilde{\mathfrak k},\tilde{H};V) \times H^{max - i}(\tilde{\mathfrak
k},\tilde{H};W) \rightarrow {\mathbb C},$$ for a dual pair of modules V
and W, and the natural duality of Zuckerman derived functors $$R^iS(V)
\times R^{max - i}S(W \otimes \varepsilon_{\tilde{\mathfrak
k}/\tilde{\mathfrak h}}) \rightarrow {\mathbb C},$$ or in a weaker form,
an isomorphism $$R^iS(V) \cong R^{max-i}S(V^\sim)^\approx.$$ It is easy
then to deduce a Borel-Weil-Bott type theorem $$R^iS
(H(L_{-w\lambda+\tilde{\rho}})) \cong
\delta_{i,s(w)}V_{-\lambda+\tilde{\rho}},$$ where $s(w) = \sum \dim
{\mathfrak g}_{\alpha} +l-1 +\sum_i a_i$ and $\tilde{\rho}$ is the sum of
the fundamental weights. Finally, it is easy to deduce an affine analogue
of the Kostant theorem on decomposition for the cohomology of the
nilpotent part $$H^i(\tilde{\mathfrak n}_+;V_{\lambda}) = \bigoplus_{w \in
\tilde{W},s(w)=i} L_{w(\lambda +\tilde{\rho})-\tilde{\rho}}$$ \vskip
1truecm

\section{Notations}

First of all we recall some notations from the theory of finite dimensional 
semisimple Lie algebras :

${\mathfrak g}$  a complex finite dimensional Lie algebra of rank $l$,

${\mathfrak h}$  some Cartan subalgebra, $dim {\mathfrak h} = l$, 

$\Delta$  the root system corresponding to the pair $({\mathfrak g},{\mathfrak h})$, 

$\Delta_+$ the system of positive roots, 

${\mathfrak n}_{\pm} = {\sum_{\alpha \in \Delta_+}} {\mathfrak g}_{\pm \alpha}$  the 
nilpotent part,

$ \Pi = \{ \alpha_1, \dots,\alpha_l\}$  the simple roots,

$\theta$  the maximal positive root,

$(.|.)$ a nondegenerate invariant bilinear form, inducing an isomorphism $ \nu :
{\mathfrak h} \rightarrow {\mathfrak h}^* $ and a form on ${\mathfrak h}^*$ normed
by the following condition $(\theta|\theta) = 2$. Hence $(\alpha|\alpha) = 2$ for
all long roots from $\Delta$,

${\check \alpha} = {2\alpha \over (\alpha|\alpha)}$ the dual root corresponding to
$\alpha$,

$\theta = \sum_{i=1}^l a_i \alpha_i$, where $a_1,\dots, a_l \in {\mathbb N}_+$,

${\check \theta} = \sum_{i=1}^l {\check a_i} {\check \alpha_i}$,where ${\check a_1}
,\dots, {\check a_l} \in {\mathbb N}_+$

$\Lambda_i\quad;\enskip(\Lambda_i,\alpha_j) = \delta_{ij}$ the fundamental weights,

$Q = \sum_{\alpha \in \Delta} {\mathbb Z}\alpha$  the root lattice in ${\mathfrak h}^*$, 

$L = \sum_{\alpha \in \Delta} {\mathbb Z}\alpha$ the sublattice of the long roots,

$\nu^{-1}(L) = {\check Q} = \sum_{\alpha \in \Delta} {\mathbb Z}{\check \alpha}$ the
co-root lattice,

$P = \sum_{i=1}^l {\mathbb Z}\Lambda_i$  the weight lattice,

$P_+ = \sum_{i=1}^l {\mathbb Z}_+\Lambda_i$ the positive chamber of fundamental weights, 

$\rho = \sum_{i=1}^l \Lambda_i$  the sum of fundamental weights,

$g = 1 + \rho({\check \theta}) = 1 + \sum_{i=1}^l {\check a}_i$  the Coxeter number,

$r_{\alpha}, \alpha \in \Delta$ the reflections defined by $$r_{\alpha}(\lambda) =
\lambda - 2(\lambda, {\check \alpha})\alpha \enskip,$$

$W = \langle r_{\alpha}, \alpha \in \Delta\rangle \subset \GL({\mathfrak
h}^*)$ the Weyl group generated by reflections, acting on $\Delta, Q, P,
etc,\dots$,.  \vskip 1truecm

Recall now the corresponding notations for the theory of affine Lie
algebras :

${\mathcal L }
 = {\mathbb C}[z,z^{-1}]$ the algebra of Laurent polynomials in variable $z \in
{\mathbb C} \setminus (0)$,

$\bar{\mathfrak g} = {\mathcal L} ({\mathfrak g}) := {\mathbb C}[z,z^{-1}]
\otimes_{\mathbb C} {\mathfrak g}$ the Laurent loop algebra,

$\langle .|. \rangle_z := I \otimes (.|.)$  the invariant Hermitian structure,

$c(X,Y) := {\Res_{z=0}} \langle{d \over {dz}}X|Y\rangle $  the
antisymmetric
2-cocycle,

$\hat{\mathfrak g}$ the central extension defined by the cocycle $c(.,.)$ and the
central element $c$, $$0 \rightarrow {\mathbb C}.c \rightarrow \hat{\mathfrak g}
\rightarrow \bar{\mathfrak g} \rightarrow 0,$$

$d = d \otimes 1 \in \Der{\mathbb C}[z,z^{-1}] \hookrightarrow \Der{\bar{\mathfrak
g}}$, acting on $\bar{\mathfrak g}$ as $z{d \over dz}$ and commuting with the
central element $c$,

$\tilde{\mathfrak g} = {\mathbb C}[z,z^{-1}] \otimes_{\mathbb C} {\mathfrak g}
\oplus {\mathbb C}c \oplus {\mathbb C}d$ the affine Lie algebra, $$[X(.) +
{\alpha}c + {\beta}d,Y(.) + {\alpha_1}c + {\beta_1}d] := ([X,Y](.) + {\beta}z{d
\over dz}Y - {\beta_1}z{d \over dz}X) + {\Res_{z=0}} \langle{d \over
dz}X|Y\rangle_zc,$$ $$\langle X(.)
+ {\alpha}c + {\beta}d|Y(.) + {\alpha_1}c + {\beta_1}d\rangle =
{\Res_{z=0}}
(z^{-1}\langle X|Y\rangle_z) + \alpha\beta_1 + \alpha_1\beta,$$

$\tilde{\mathfrak h} = {\mathfrak h} \oplus {\mathbb C}c \oplus {\mathbb C}d$ some
Cartan affine subalgebra, where ${\mathfrak h}$ is some ( finite dimensional)
Cartan subalgebra of ${\mathfrak g}$,

$\tilde{\mathfrak g}_{\tilde{\mathfrak \alpha}} = \{X \in \tilde{\mathfrak g} ;
[H,X] = \tilde{\alpha}(H)X, H \in \tilde{\mathfrak h} \}$ a root space,

$\tilde{\Delta} = \{ \tilde{\alpha} ; \tilde{\mathfrak g}_{\tilde{\alpha}} \ne 0\} =
\tilde{\Delta}_W \cup \tilde{\Delta}_I$ the root system,

$\tilde{\Delta}_W = \{ \tilde{\alpha} = \alpha + n\delta ; \alpha \in \Delta, n \in
{\mathbb Z} \}$ the real roots,

$\delta$ such a root that $\delta(d) = c^*(d) = 1, \delta|_{{\mathfrak h} +
{\mathbb C}c} = 0$,

$\tilde{\Delta}_I = \{ n\delta ; n \in {\mathbb Z} \setminus (0) \}$ the imaginary
roots, $$\dim_{\mathbb C} \tilde{\mathfrak g}_{\tilde{\alpha}} = \cases{1 &\mbox{ if
}$\alpha \in \tilde{\Delta}_W$ \cr l &\mbox{ if } $\alpha \in \tilde{\Delta}_I,$
\cr}$$

$\tilde{\mathfrak n}_{\pm} = {\mathfrak n}_{\pm} \oplus \sum_{n>0}^{\oplus} z^n
\otimes {\mathfrak g}$ the nilpotent parts,

$\tilde{\mathfrak g} = \tilde{\mathfrak h} \oplus \tilde{\mathfrak n}_+ \oplus
\tilde{\mathfrak n}_- = \tilde{\mathfrak h} \oplus \sum^{\oplus}_{\tilde{\alpha} \in
\tilde{\Delta}} \tilde{\mathfrak g}_{\tilde{\alpha}}$ the Cartan decomposition,

$\delta = d^* \in \tilde{\mathfrak h}^*, \delta(d) = c^*(d) = 1,
\delta|_{{\mathfrak h} \oplus {\mathbb C}c} =0$,

$\tilde{\Pi} = \{ \alpha_0 := \delta - \theta, \alpha_1,\dots \alpha_l \}$ the
simple roots,

$\check{\tilde{\Pi}} = \{ \check{\alpha}_0 := c -
\check{\theta},\check{\alpha}_1,\dots, \check{\alpha}_l \}$ the simple co-roots,

$\{ \tilde{\Lambda}_0, \dots, \tilde{\Lambda}_l ;\tilde{\Lambda}_i = \Lambda_i +
{\check a}_i\Lambda_0 \}$ the fundamental weights, normed as
$\langle\tilde{\Lambda}_i |
{\check \alpha}_j\rangle = \delta_{ij}$,

$\tilde{\Delta}_+ = \Delta_+ \cup \{ k\delta +\alpha ; \alpha \in \Delta, k \in
{\mathbb Z}_+ \} \cup \{ k\delta ; k \in {\mathbb Z}_+ \}$,

$\tilde{P} = \sum_{i=0}^l {\mathbb Z}\tilde{\Lambda}_i$  the weight lattice, 

$\tilde{P}_+ = \sum_{i=0}^l {\mathbb Z}_+ \tilde{\Lambda}_i$ the chamber of the dominant weights,

$\tilde{\rho} = \sum_{i=0}^l \tilde{\Lambda}_i = \rho + g\tilde{\Lambda}_0$  the sum of the fundamental weights,

$\tilde{Q} = {\mathbb Z}\delta \oplus Q$  the affine root lattice,

$\tilde{A} = ( \alpha_j({\check \alpha_i}))_{i,j = 0,1,\dots,l}$  the generalized Cartan matrix, 

$\tilde{W} = W \ltimes {\check Q}$  the affine Weyl group, generated by the reflections $r_{\tilde{\alpha}},\tilde{\alpha} \in \tilde{\Delta}$,

$V_{\lambda}$ the standard $\tilde{\mathfrak g}$-module with highest
weight $\lambda \in \tilde{P}$ and with the weight vector $v_{\lambda} \in
V_{\lambda}$ such that $$\tilde{\mathfrak h}.v_{\lambda} =
\lambda(\tilde{\mathfrak h}).v_{\lambda},$$ 
$$\tilde{\mathfrak n}_+.v_{\lambda} = 0 $$
and such that for all $v \in V_{\lambda}, \dim U(\tilde{\mathfrak n}_{\pm})v < \infty$,

$U(\tilde{\mathfrak g})$  the universal enveloping algebra and 

${\mathcal Z}
(\tilde{\mathfrak g})$  the center of $U(\tilde{\mathfrak g})$. 
\vskip 1truecm 

\section{Category ${\mathcal C}(\tilde{\mathfrak k},\tilde{H})$}

Assume $\tilde{\mathfrak g} = \tilde{\mathfrak k}$ to be an affine Lie algebra associated to a compact complex Lie algebra ${\mathfrak k = k}_{\mathbb R} \otimes_{\mathbb R} {\mathbb C}$, where ${\mathfrak k}_{\mathbb R}$ is the real Lie algebra of some 
compact connected Lie group
$K$. Suppose that 2-cocycle $c(.,.)$ is integral. In this case the associated
Lie algebra extension can be lifted to an extension of the corresponding Lie groups
$$1 \rightarrow {\mathbb T} \rightarrow \tilde{K} \rightarrow {\mathcal L} K \cong
C^{\infty}({\mathbb S}^1, K) \rightarrow 1.$$ For each subgroup $H \subset K$,
consider the corresponding central extension $$1 \rightarrow {\mathbb T} \rightarrow
\tilde{H} \rightarrow H \rightarrow 1.$$ In particular, if $H$ is some Cartan
subgroup of $K$ we have an affine extension $\tilde{H}$ of Cartan subgroup $H$. In
this case, the affine Cartan subgroup $\tilde{H}$ is also finite dimensional.

The pair $(\tilde{\mathfrak k},\tilde{H})$ is called to be {\it compatible, } if
there is such a homomorphism of affine Lie groups $\tilde{H} \rightarrow
\Aut\tilde{\mathfrak k}$ that the corresponding homomorphism of affine Lie algebras
$\tilde{\mathfrak h} \rightarrow \End\tilde{\mathfrak k}$ is coincided with the
adjoint representation $\ad_{\tilde{\mathfrak h}}(\tilde{\mathfrak k})$ and that it
induces the adjoint representation $\Ad_{\tilde{\mathfrak h}}(\tilde{H})$ in
$\tilde{\mathfrak k}$. 

For such a fixed compatible pair $(\tilde{\mathfrak k},\tilde{H})$ we define
${\mathcal C} (\tilde{\mathfrak k},\tilde{H})$ as the category the objects of which
have simultaneously $(\tilde{\mathfrak k},\tilde{\mathfrak h})$-module (see [3]) and
$\tilde{H}$-module structures, which are compatible in the sense $$h(Xv) = ((\Ad
h)X).hv,\forall h \in \tilde{H}, \forall v \in V, X \in \tilde{\mathfrak k}.$$

Remark that all the categories ${\mathcal C} (\tilde{\mathfrak h},\tilde{H}),
{\mathcal C} (\tilde{\mathfrak k},\tilde{H}),{\mathcal C} (\tilde{\mathfrak
k},\tilde{K}), \dots$, raised in the introduction are the particular cases of this
definition. Following M. Duflo and M. Vergne [2] we define now the derived  Zuckermann
functors and then prove their duality in the category ${\mathcal C }
(\tilde{\mathfrak k}, \tilde{H}).$

  Denote by $F(\tilde{K}) = \oplus_{\lambda \in \tilde{P}} (V^*_{\lambda} \otimes
V_{\lambda})$ the associated algebra of matrix elements of the standard modules,
where the sum runs over the set of weight lattice $P$. It is also the
$(\tilde{\mathfrak k },\tilde{H})$-module of regular functions on $\tilde{K}$. 
Denote by $r$ and $l$ the right and left regular representations of $\tilde{K}$ ( or
$\tilde{\mathfrak k}$ ) on $F(\tilde{K})$, respectively. Let $V$ be a
$(\tilde{\mathfrak k},\tilde{H})$-module with the action $\check{\theta}$ of
$\tilde{\mathfrak k}$ ( or $\tilde{H} )$ in $V$. Then we define a $(\tilde{\mathfrak
k},\tilde{H})$-module structure on $F(\tilde{K};V) := F(\tilde{K}) \otimes V$ by
$$(\psi(X)f)(k) = \check{\theta}(X) f(k) + (l(X)f)(k), \forall x \in
\tilde{\mathfrak k}, \forall k \in \tilde{K},$$ $$(\psi(h)f)(k) =
\check{\theta}(h) f(h^{-1}k), \forall k,h \in \tilde{H}.$$

This action commutes with the $(\tilde{\mathfrak k},\tilde{H})$-module structure
$\psi$. The right regular representation $r$ ( of algebra $\tilde{\mathfrak k}$, or
of group $\tilde{K}$ ) also commutes with the $(\tilde{\mathfrak
k},\tilde{H})$-module  structure
$\psi$ and finally $$r(k)\psi(X)r(k^{-1}) = \Gamma(k.X).$$ The representation
$\Gamma$ of the Algebra $\tilde{\mathfrak k}$ and the representation $r$ of the
group $\tilde{K}$ can be continued on $$\Hom(\wedge^* \tilde{\mathfrak k},
F(\tilde{K};V)) \cong F(\tilde{K};\Hom(\wedge^* \tilde{\mathfrak k},V)).$$
 
It is easy to prove ( see [2], p. 468 ) that for every $Y \in \tilde{\mathfrak k}$,
one can define $I(Y)$ by $$ (I(Y)\mu)(k) = i(k.Y)\mu(k), \forall k \in \tilde{K},
\forall \mu \in F(\tilde{K}, \Hom (\wedge^* \tilde{\mathfrak k},V)) $$ and that
$\Gamma(Y) - r(Y) = d{\circ}I(Y) + I(Y){\circ}d$. Thus $\Gamma$ and $r$ induces a
$(\tilde{\mathfrak k},\tilde{K})$-module structure on the Zuckermann derived functors
$$\Gamma^i(V) := \bigoplus_{\lambda \in \tilde{P}} H^i(\tilde{\mathfrak
k},\tilde{H};  V^*_{\lambda} \otimes V) \otimes V_{\lambda}.$$ In the next
section, by using of the general ideas of homological algebra [3], we show that for
connected $\tilde{K}, \Gamma^{\circ}(V)$ coincides with the Zuckermann functor $S(V)
:= V[\tilde{\mathfrak k}]$ from the category ${\mathcal C} (\tilde{\mathfrak
k},\tilde{H})$ into the category ${\mathcal C} (\tilde{\mathfrak k},\tilde{K})$. 

Remark that in the second part of this contribution we shall prove a version of the
Peter-Weyl theorem, which asserts that $$ S (U(\tilde{\mathfrak k})^*) \cong
\bigoplus_{\lambda \in \tilde{P}} V^*_{\lambda} \otimes V_{\lambda} =
F(\tilde{K}).$$ Thus $F(\tilde{K})$ can be interpreted as the algebra of regular
functions on $\tilde{K}$, as in the classical case, \cite{kp}, \cite{t}. 

Recall that $\tilde{\mathfrak k}/\tilde{\mathfrak h}$ is infinite
dimensional space ${\mathcal L} [z,z^{-1}] \otimes ({\mathfrak
k}/{\mathfrak h})$. But it is well-defined the determinant fiber bundle
structure on $ \det (\tilde{\mathfrak k}/\tilde{\mathfrak h}) :=
\wedge^{max} (\tilde{\mathfrak k}/\tilde{\mathfrak h}) $ ( see for example
[5] ), in case where ${\mathfrak h}$ is a Cartan subalgebra. Because the
corresponding group $H$ is compact, $\det (\tilde{\mathfrak
k}/\tilde{\mathfrak h})$ admits also a $(\tilde{\mathfrak
 k},\tilde{H})$-module structure $\varepsilon_{\tilde{\mathfrak k}/\tilde{\mathfrak
h}}$ with the trivial differential.

Suppose that $\langle .,.\rangle$ is a $(\tilde{\mathfrak k},\tilde{H})$-invariant
pairing of 
$(\tilde{\mathfrak k},\tilde{H})$-modules $V,W \in {\mathcal C}
(\tilde{\mathfrak k},\tilde{H})$. 

\begin{theorem} For each $(\tilde{\mathfrak k},\tilde{H})$-invariant pairing of $V,W
\in {\mathcal C} (\tilde{\mathfrak k},\tilde{H})$, there exists a $(\tilde{\mathfrak
k},\tilde{K})$-invariant pairing of type $$\langle \Gamma^i(V),\Gamma^{\max - i}(W
\otimes
\varepsilon_{\tilde{\mathfrak k},\tilde{\mathfrak h}})\rangle \rightarrow {\mathbb
C}.$$
\end{theorem}
{\it Proof}.  Let us denote by $dW_K$ the Wiener measure on $\tilde{K}$. 
Then each $\tilde{H}$-invariant $i^{th}$ relative cocycle $$\omega_1 \otimes f_1 \in
\wedge^i (\tilde{\mathfrak k}/\tilde{\mathfrak h})^* \otimes F(\tilde{K};V)$$ can be
$(\tilde{\mathfrak k},\tilde{H})$-invariantly paired with any relative cocycle
$\omega_2 \otimes f_2 \otimes \zeta$ of degree 1 from $$\wedge^{max -
i}(\tilde{\mathfrak k}/\tilde{\mathfrak h})^* \otimes F(\tilde{K};W) \otimes
\wedge^{max}(\tilde{\mathfrak k},\tilde{\mathfrak h})$$ by the formula
$$\langle\omega_1
\otimes f_1, \omega \otimes f_2 \otimes \zeta\rangle = \langle\omega_1 \wedge
\omega_2,\zeta\rangle
\int_{\tilde{K}} \langle f_1(k),f_2(k)\rangle  dW_K(k).$$ This pairing is
compatible with
the differential $d$ of the complex of relative $(\tilde{\mathfrak
k},\tilde{H})$-valued cocycles, i.e. for $\tilde{H}$-invariant elements $\mu_1 \in
\wedge^i(\tilde{\mathfrak k}/\tilde{\mathfrak h})^* \otimes F(\tilde{K};  V)$ and
$\mu_2 \in \wedge^{max - i}(\tilde{\mathfrak k},\tilde{\mathfrak h})^* \otimes
F(\tilde{K}) \otimes \wedge^{max}(\tilde{\mathfrak k},\tilde{\mathfrak h})$, we have
$$\langle d\mu_1,\mu_2\rangle + (-1)^i\langle\mu_1,d\mu_2\rangle = 0.$$ The last
gives thus a
$(\tilde{\mathfrak k},\tilde{K})$-invariant pairing of the Zuckermann derived
functors of type $$\langle \Gamma^i(V), \Gamma^{max - i}(W \otimes
\varepsilon_{\tilde{\mathfrak k},\tilde{\mathfrak h}})\rangle \rightarrow {\mathbb C}.$$ 
\hfill$\square$

\section{Duality of the Zuckermann derived functors}

We develop in this section the formal homological algebra on the categories ${\mathcal C}
(\tilde{\mathfrak h},\tilde{H})$, ${\mathcal C}(\tilde{\mathfrak k},\tilde{H}),
{\mathcal C}(\tilde{\mathfrak k},\tilde{K})$ etc.

First of all, it is easy to see that {\it the category ${\mathcal C}
(\tilde{\mathfrak k},\tilde{H})$ has enough injective objects. } Really, for every
object $W$ of ${\mathcal C} (\tilde{\mathfrak k},\tilde{H})$, the object $I(W) =
\Hom_{U(\tilde{\mathfrak h})}(U(\tilde{\mathfrak k}),W)$ is injective in the
category
${\mathcal C} (\tilde{\mathfrak k}, \tilde{H})$. Furthermore every injective object
in ${\mathcal C} (\tilde{\mathfrak k},\tilde{H})$ can be presented as a direct
summand of an object of this type $I(W)$. More exactly, we have $V \hookrightarrow
I(V)$. It is easy to see that if $F \in {\mathcal C} (\tilde{\mathfrak
k},\tilde{K})$, then $F \in {\mathcal C} (\tilde{\mathfrak k},\tilde{H})$ and the
exact covariant additive functor $F \otimes (.)$ maps the injective objects into the
same ones, commuting with the Zuckermann functor $S(V) = V[\tilde{\mathfrak k}]$. 
Hence, we have the isomorphism $$R^iS (F \otimes V) \cong F \otimes R^iS (V).$$ In
particular, to each resolution $0 \rightarrow V \rightarrow I^*$ of type $$ 0
\rightarrow V \rightarrow I^0 \rightarrow I^1 \rightarrow \dots $$ we can apply the
functor $U(\tilde{\mathfrak k}) \otimes (.)$ and obtain the commutativity of $R^iS$
with $U(\tilde{\mathfrak k}) \otimes (.)$ and also with the multiplication functor
$m$ (i.e. the $U(\tilde{\mathfrak k})$-action ), following the commutative diagram
$$\begin{picture}(100,120)(10,20)
\put(-110,110){$U(\tilde{\mathfrak k})\otimes R^iS(V)$}
\put(150,110){$R^iS(U(\tilde{\mathfrak k}) \otimes V)$}
\put(20,20){$R^iS(V)$}
\put(-15,110){\vector(1,0){150}}
\put(40,115){$\cong$}
\put(-75,90){\vector(3,-2){80}}
\put(-45,55){m}
\put(70,40){\vector(3,2){80}}
\put(110,55){$R^iS(m)$}
\end{picture}$$
More particularly, if the submodule $V \subset V_1
\oplus V_2 \in {\mathcal C} (\tilde{\mathfrak k},\tilde{H})$ is stable under
${\mathcal Z} (\tilde{\mathfrak k}) = U(\tilde{\mathfrak k})^{\tilde{\mathfrak k}}$
and if for some element $ u \in {\mathcal Z} (\tilde{\mathfrak k})$, the operator $u
- \chi(u)$ acts nilpotently on $V$, then so it acts also on the derived functors
$R^iS(V)$. 

Consider the finite dimensional injective $\tilde{H}$-invariant Koszul resolution
$$0 \rightarrow {\mathbb C} \rightarrow \Hom_{U(\tilde{\mathfrak
h})}(U(\tilde{\mathfrak k}), \wedge^0(\tilde{\mathfrak k}/\tilde{\mathfrak h})^*)
{\buildrel d \over \rightarrow} \dots $$ $$\dots {\buildrel d \over \rightarrow}
\Hom_{U(\tilde{\mathfrak h})}(U(\tilde{\mathfrak k}),\wedge^{max}(\tilde{\mathfrak
k}/\tilde{\mathfrak h})^*) \rightarrow 0$$ and apply the tensor multiplication with
$V$ on the right. We have then an injective resolution of $V$ $$0 \rightarrow V
\rightarrow \Hom_{U(\tilde{\mathfrak h})}(U(\tilde{\mathfrak k}),
\wedge^0(\tilde{\mathfrak k}/\tilde{\mathfrak h})^* \otimes V) {\buildrel d_V \over
\rightarrow} \dots $$ $$\dots {\buildrel d_V \over \rightarrow }
\Hom_{U(\tilde{\mathfrak h})}(U(\tilde{\mathfrak k}),\wedge^{max}(\tilde{\mathfrak
k}/\tilde{\mathfrak h})^* \otimes V) \rightarrow 0.$$ Applied the Zuckermann functor
$S(V) = V[\tilde{\mathfrak k}]$ and taken the $\tilde{H}$-invariant parts, we have
the derived functors as the cohomology groups of complex $$R^iS(V) =
\Coh^i(I(\wedge^*(\tilde{\mathfrak k}/\tilde{\mathfrak h})^* \otimes V)^{
\tilde{H}},d_V).$$ Remark that $$\Hom_{U(\tilde{\mathfrak h})}(U(\tilde{\mathfrak
k}),\wedge^i(\tilde{\mathfrak k}/\tilde{\mathfrak h}) \otimes V) \cong
\Hom_{U(\tilde{\mathfrak h})} (\wedge^i(\tilde{\mathfrak k}/\tilde{\mathfrak
h}),U(\tilde{\mathfrak k})^* \otimes V).$$ Hence, applying the functor $S$, we have
$$I(\wedge^i(\tilde{\mathfrak k}/\tilde{\mathfrak h})^* \otimes V) \cong
\Hom_{U(\tilde{\mathfrak h})}(\wedge^i(\tilde{\mathfrak k}/\tilde{\mathfrak
h}),R(\tilde{\mathfrak k}) \otimes V),$$ where $R(\tilde{\mathfrak k}) =
\bigoplus_{\lambda \in \tilde{P}} V_{\lambda} \otimes V^*_{\lambda}$ is the maximal
$(\tilde{\mathfrak k},\tilde{H})$-bimodule in $U(\tilde{\mathfrak k})^*$, i.e.
$$I(\wedge^i(\tilde{\mathfrak k}/\tilde{\mathfrak h})^* \otimes V) \cong
\bigoplus_{\lambda \in \tilde{P}} \Hom_{U(\tilde{\mathfrak
h})}(\wedge^i(\tilde{\mathfrak k}/\tilde{\mathfrak h}),V \otimes V^*_{\lambda})
\otimes V_{\lambda}.$$ Therefore, we have the equivalence of functors $$R^iS(V)
\cong H(V) := \bigoplus_{\lambda \in \tilde{P}} H^i(\tilde{\mathfrak
k},\tilde{H};V_{\lambda} \otimes V) \otimes V^*_{\lambda},$$ $$R^{max - i}S(V)
\cong G(V) := \bigoplus_{\lambda \in \tilde{P}} H^{max - i}(\tilde{\mathfrak
k},\tilde{H};V^\sim \otimes V_{\lambda}) \otimes V^*_{\lambda}.$$ Applying
Theorem 1, we have now an equivalence of functors $$T_V : R^iS(V) \quad {\buildrel
\simeq \over \rightarrow} \quad R^{max - i}S(V^\sim)^\sim.$$ Hence we obtain the
result

\begin{theorem} {\it $V \longrightarrow T_V$ is an  equivalence of
functors.}
\end{theorem}

\section{Category ${\mathcal O}$} 

By using the trivial action of the nilpotent part $\tilde{\mathfrak n}_+$ of affine
Lie algebras $\tilde{\mathfrak b} := \tilde{\mathfrak h} \oplus \tilde{\mathfrak
n}_+$, one can, in place of the category ${\mathcal C} (\tilde{\mathfrak
k},\tilde{H})$ and the functor $I$ consider the category ${\mathcal C}
(\tilde{\mathfrak b},\tilde{H})$ and the functor $H$, $$H(W) :=
\Hom_{U(\tilde{\mathfrak b})}(U(\tilde{\mathfrak k}),W).$$ Thus we can use the
subcategory ${\mathcal O}$ in ${\mathcal C} (\tilde{\mathfrak k},\tilde{H})$, which
has also enough injective objects. This is a {\it more economic } way of computing
the derived functors. Remark that the most important Verma modules belong too to
this subcategory. We have deal with this situation in proving the Borel-Weil-Bott
type theorem. 

For each $H(W) \in {\mathcal O}$, where $W \in {\mathcal C} (\tilde{\mathfrak
b},\tilde{H})$, let us consider the ${\mathcal O}$-injective, $\tilde{H}$-invariant
resolution of type $$0 \rightarrow H(W) \rightarrow I^*(W),$$ Lie algebra
$$\tilde{\mathfrak n} = \tilde{\mathfrak n}_+ \oplus \tilde{\mathfrak n}_-, $$
and $$I^i(W) := (\Hom_{U(\tilde{\mathfrak h})}(U(\tilde{\mathfrak
k}),\wedge^i{\mathfrak n}^* \otimes W)[\tilde{\mathfrak h}])^{\tilde{H}} \in
{\mathcal O}.$$ Then we have $$SI^i(W) \cong \Hom_{U(\tilde{\mathfrak
h})}(\wedge^i\tilde{\mathfrak n},R(\tilde{\mathfrak k}) \otimes
W)^{\tilde{\mathfrak{H}}} $$
  $$\cong \bigoplus_{\lambda \in \tilde{P}} V_{\lambda} \otimes (
\Hom_{U(\tilde{\mathfrak h})} ( \wedge^i\tilde{\mathfrak n}, V_{\lambda}^* \otimes
W))^{\tilde{H}}.$$ Hence $$R^i S(H(W)) \cong \bigoplus_{\lambda \in \tilde{P}}
V_{\lambda} \otimes (H^i(\tilde{\mathfrak n};V^*_{\lambda} \otimes W))^{\tilde{H}}
.$$ In the category ${\mathcal O }$ an analogue of the classical
Bernstein-Gel'fand-Gel'fand's results [1] holds :{\it Every indecomposable injective
object is a direct summand of an object of type } $$H((U(\tilde{\mathfrak
n})/\tilde{\mathfrak n}^r U(\tilde{\mathfrak n}))^* \otimes L_{\lambda}), \lambda
\in \tilde{P}, \forall r >> 0.  $$ Therefore, remarking that if $F \in {\mathcal C}
(\tilde{\mathfrak k},\tilde{H})$, $$H^i(\tilde{\mathfrak n},(U(\tilde{\mathfrak
n})/\tilde{\mathfrak n}^r U(\tilde{\mathfrak n}))^* \otimes L_{\lambda} \otimes F) =
0, i = 1,2,...$$ we see that in the category ${\mathcal O}$ all the objects are
acyclic, i.e. if $V$ is an injective object in ${\mathcal O}$, then $R^iS(V)$
vanishes unless $i = 0 $. In the category ${\mathcal O}$, all the objects have
finite injective cohomological dimension. Then by induction on this dimension one
proves that $(R^iS)|_{{\mathcal O} } \cong R^i(S|_{\mathcal O}).$ This mean that
computing the Zuckermann derived functors for objects in ${\mathcal O}$ is
independent from the bigger category ${\mathcal C} (\tilde{\mathfrak k},\tilde{H})$.

\section{Borel-Weil-Bott and Kostant type theorems}

Recall that by $\tilde{P}$ we denote the weight lattice of our compact complex
affine Lie algebra, $\tilde{P}_+$ the Weyl chamber of the dominant weights. The Weyl
group $\tilde{W}$ acts on $\tilde{P}$ and every regular weight $\mu$ can be uniquely
presented in form $-w\lambda + \tilde{\rho}$, where $\lambda \in \tilde{P}_+$ and $w
\in \tilde{W}$. Recall that the length function is defined by $l(w) := \# \{ \lambda
\in \tilde{\Pi} ; w\lambda = -\lambda \} \leq l+1$. If $\varpi \in W$ and $l(\varpi)
= l + 1$, then $ \varpi\tilde{\Pi} = -\tilde{\Pi}$. 

By definition, we have $\tilde{\rho} = \sum_{i=0}^l \lambda_i$ and
$$-\varpi\tilde{\rho} + \tilde{\rho} = 2\sum_{i = 0}^l \lambda_i.$$ Remark
that although $\dim_{\mathbb C} {\mathfrak g}_{k\delta} = l, \forall k \in {\mathbb
Z} \setminus (0)$, where $\delta\in {\mathfrak h}^*, \delta(d) = 1$ and
$\delta|_{{\mathfrak h}+ {\mathbb C}c} = 0,$ we have always $\dim
L_{-\varpi\tilde{\rho} + \tilde{\rho}} =1$. 

Let us introduce a function 
$$s(w) = \left\{ \begin{array}{ll} l(w) &\mbox{ if } w\alpha_0  \ne -\alpha_0 \\
                l(w) +l - 1 + \sum_i a_i &\mbox{ if } w\alpha_0 = -\alpha_0.
\end{array}\right. $$

\begin{theorem}\label{Thm6.1}
For every regular integral dominant weight $\lambda \in \tilde{P}_+ $
and for every element $w \in \tilde{W}$,
$$R^iS(H(L_{-w\lambda + \tilde{\rho}})) = \cases{0 &\mbox{ if } $i \ne s(w)$ \cr
  V_{-w\lambda + \tilde{\rho}} &\mbox{ if } $i = s(w)$. \cr}$$ \end{theorem}
{\it Proof}. As in the classical case, we prove the theorem by induction
on the length of elements of the affine Weyl group $\tilde{W}$. 

(a) First of all we verify the assertion in the case of maximal length elements,
$l(w) = l+ 1, w\tilde{\Pi} = -\Pi $. Hence $w = \varpi$.

From the discussed properties of the functor $H$ one deduces easily that
$$H(L_{-\varpi\lambda}) \cong H(L_{-\varpi\tilde{\rho} + \tilde{\rho}}
\otimes V_{-\lambda + \tilde{\rho}})_{\chi_{\lambda}},$$ where the index
$\chi_{\lambda}$ means the submodule, where the operators $z -
\chi_{\lambda}(z)$ act nilpotently, for all $z \in {\mathcal Z}
(\tilde{\mathfrak k}) = \cent U(\tilde{\mathfrak k})$, and
$\chi_{\lambda}$ the infinitesimal character corresponding to $\lambda$. 

In Section 3 we have seen that 
$$R^iS(W_{\chi_{\lambda}}) \cong  (R^iS(W))_{\chi_{\lambda}}, W \in {\mathcal C}
(\tilde{\mathfrak k},\tilde{H}),$$
$$R^iS(W \otimes V_{\lambda}) \cong R^iS(W) \otimes V_{\lambda}. $$
Hence, we have 
$$R^iS(H(L_{-\varpi\lambda + \tilde{\rho}})) \cong
R^iS((H(L_{-\varpi\tilde{\rho} + \tilde{\rho}}) \otimes V_{-\lambda
+\tilde{\rho}}))_{\chi_{\lambda}}$$
 $$\cong (R^iS(H(L_{-\varpi\tilde{\rho} + \tilde{\rho}}) \otimes
V_{-\lambda + \tilde{\rho}}))_{\chi_{\lambda}} \cong
(R^iS(H(L_{-\varpi\tilde{\rho} + \tilde{\rho}})) \otimes V_{-\lambda +
\tilde{\rho}})_{\chi_{\lambda}}. $$ Therefore it is enough to compute
$R^iS(H(L_{-\varpi\tilde{\rho} + \tilde{\rho}}))_{\chi_{\lambda}}$. We
have $$R^iS(H(L_{\varpi\tilde{\rho}+\tilde{\rho}})) = (\bigoplus_{\lambda
\in \tilde{P}} V_{\lambda} \otimes H^i(\tilde{\mathfrak
n};L_{-\varpi\tilde{\rho}+\tilde{\rho}} \otimes V^*))_{\chi_0}$$
  $$= V_0 \otimes H^i(\tilde{\mathfrak
n};L_{-\varpi\tilde{\rho}+\tilde{\rho}} \otimes V^*_0)^{\tilde{H}} 
 \cong V_0 \otimes H^i(\tilde{\mathfrak
n};L_{-\varpi\tilde{\rho}+\tilde{\rho}}),$$ where $\tilde{\mathfrak n}
:= \tilde{\mathfrak n}_+ \oplus \tilde{\mathfrak n}_-$,
$\wedge^i{\mathfrak n}_+$ has no vectors of weight
$-\varpi\tilde{\rho}+\tilde{\rho}$, if $i \ne s(w)$. If $i = s(w)$, we
have $$\wedge^{s(w)}\tilde{\mathfrak n} = (\wedge^{s(w)}\tilde{\mathfrak
n})_{-\varpi\tilde{\rho}+\tilde{\rho}} \oplus \dots$$ and
$$H^{s(w)}(\tilde{\mathfrak n};L_{-\varpi\tilde{\rho}+\tilde{\rho}}) =
{\mathbb C}.$$ The theorem is verified in this case. Remark that here we
have the non-triviality of the cohomology group in dimension $s(w)$, while
in the classical case the non-triviality in the dimension $l(w)$. 

(b) Suppose that $l(w) < l+1$ and that the theorem is verified for all $w' \in
\tilde{W}$, such that $$l(w') \geq l(w) + 1.$$ Consider $\tilde{Q} =
w\tilde{P}$. Then $w\lambda$ is $\tilde{Q}$-dominant, integral and regular. Because
$\tilde{Q} \ne w\tilde{\Delta}_+$, there exists a $\tilde{Q}$-simple root $\alpha$
from $\tilde{\Delta}_+$, such that $\langle\lambda,\alpha\rangle > 0$. Because
$$s_{\alpha}\tilde{Q} \supseteq \tilde{Q} - \{\alpha\},$$ then $w\lambda$ and
$s_{\alpha}w\lambda$ are $P_1$-dominant, where $\tilde{P}_1 = \tilde{Q} -
\{\alpha\}$.

\begin{lemma} Suppose that $\lambda \in \tilde{P}$ is regular integral weight. Then
there exists a unique root system $Q$ such that $\lambda$ is dominant weight. 
Suppose that ther exists a $Q$-simple root $\alpha \in \tilde{\Delta}_+ \cup Q$ such
that $\langle\lambda, \alpha\rangle > 0$ and that the both elements $\lambda$ and $\lambda' =
s_{\alpha}\lambda$ are $(Q-\{\alpha\})$-dominant. Then in the category ${\mathcal O
}$ there exists such an object $E$ that :

(i) the sequence $$0 \rightarrow M_{\lambda} \rightarrow E \rightarrow M_{\lambda'}
\rightarrow 0$$ is exact and

(ii) $R^iS(E) = 0, \forall i = 0,1,2,\dots $, where $M_{\lambda} =
U(\tilde{\mathfrak k}) \otimes_{U(\tilde{\mathfrak b})} L_{\lambda - \tilde{\rho}}$
is the Verma modules associated to $\lambda$. 
\end{lemma}
{\it Proof}.
This lemma is an affine analogue of Lemma 6.2 from [3] and its proof does not
require an essential change.\hfill$\square$

\vskip 1cm

{\sl End of the proof of Theorem \ref{Thm6.1}}. 

From this lemma, by changing $\lambda \mapsto w\lambda$, there exists a module $E
\in {\mathcal O} $ satisfying (i) and (ii). Then $E$ has the finite dimensional
weight spaces with respect to $\tilde{\mathfrak h}$.  Hence we have an exact
sequence $$ 0 \rightarrow M^{\sim}_{s_{\alpha}w\lambda} \rightarrow E^\sim
\rightarrow M^{\sim}_{w\lambda} \rightarrow 0.$$ In virtue of the duality theorem,
we have $$R^iS(E^{\sim}) = R^iS(E)^{\approx} = 0,\forall i \in {\mathbb N}.$$ From
the long exact sequence, one has $$R^iS(M^{\sim}_{w\lambda}) \cong
R^{i+1}S(M^{\sim}_{s_{\alpha}w\lambda}).$$ Remark that $$M^{\sim}_{w\lambda} =
H(L_{-w\lambda + \tilde{\rho}}), $$ $$M^{\sim}_{s_{\alpha}w\lambda} =
H(L_{-s_{\alpha}w\lambda + \tilde{\rho}}) \quad.$$ Hence $$R^iS(H(L_{-w\lambda +
\tilde{\rho}})) \cong R^{i+1}S(H(L_{-s_{\alpha}w\lambda + \tilde{\rho}})) = \left\{
\begin{array}{ll} 0 &\mbox{ if } i + 1 \ne s(s_{\alpha}w) \\
 V_{-\lambda + \tilde{\rho}} &\mbox{ if } i + 1 =
s(s_{\alpha}w),\end{array}\right.$$
i.e. $$R^iS(H(L_{-w\lambda + \tilde{\rho}})) \cong \left\{\begin{array}{ll} 0
&\mbox{ if } i \ne s(w)\\ V_{-\lambda + \tilde{\rho}} &\mbox{ if } i = s(w).
\end{array}\right.$$ The theorem is completely proved. \hfill $\square$

As usually, changing $\lambda$ by $\lambda + \tilde{\rho}$, we have 
$$\delta_{i,s(w)}V_{-\lambda} = R^iS(H(L_{-w(\lambda + \tilde{\rho}) +
\tilde{\rho}})) $$
 $$\cong V_{-\lambda} \otimes H^i(\tilde{\mathfrak n};L_{-w(\lambda +
\tilde{\rho}) + \tilde{\rho}} \otimes V^*_{\lambda})^{\tilde{H}} $$

   $$\cong V_{-\lambda} \otimes \Hom_{U(\tilde{\mathfrak
h})}(L^*_{-w(\lambda +\tilde{\rho}) - \tilde{\rho}}, H^i(\tilde{\mathfrak
n};V^*_{-\lambda}))^{\tilde{H}}.$$ Hence we have just {\sl an affine
analogue of the Kostant theorem} on cohomology of the nilpotent part. 

\begin{theorem} $$H^i(\tilde{\mathfrak n}_+;V) \cong \bigoplus_{w \in
\tilde{W},s(w)=i} L_{w(\lambda + \tilde{\rho}) - \tilde{\rho}}$$
\end{theorem}

\noindent{\it Remark}. We have deal with $(\tilde{\mathfrak
k},\tilde{H})$-modules and the central characters of type $e^{\lambda}$
and the infinitesimal characters of type $\chi_{\lambda}$. A similar
situation with $(\tilde{\mathfrak g},\tilde{K})$-modules with the central
characters of Harish-Chandra type $\theta_{\lambda}$ and the
infinitesimal characters $\chi_{\lambda}$ gives us an algebraic
realization of the discrete series representations for loop groups. This
idea will be devoted the next part of our contribution.

\section*{Acknowledgments}

This work is completed under the awarded to the author research fellowship
grant of Alexander von Humboldt Foundation and the perfect condition of
working at the Department of Mathematics and Sonderforschungbereich 343,
Bielefeld University, to which the author is very happy to thank. The deep
thanks are addressed to Prof.  Dr. Anthony Bak for his constant support
and consideration. The author thanks Professors M. Duflo and M. Vergne for
having brought to his knowledge their results [2] and Professor V. Kac for
reference [4].

\bibstyle{plain}

\vskip 2cm
Bielefeld, February 1992

\end{document}